\documentclass[11]{article}
\usepackage{mathrsfs}
\usepackage{latexsym,amsfonts,amssymb}
\usepackage{amsmath}
\usepackage{amsbsy}
\usepackage{epsfig}
\usepackage{enumerate}

\newtheorem{theorem}{Theorem}

\newtheorem{remark}{Remark}

\newbox\TempBox \newbox\TempBoxA

\def\ep{\textsf{E}} 
\def\Var{\textsf{Var}} 
\def\Cal#1{{\mathcal #1}}

\begin{document}

\title{\huge \bf A note on the invariance principle of the product of sums of random variables
\footnote {Research supported by Natural Science Foundation of China
(NSFC) (No. 10471126).
\newline  $^{1,2}$Department of Mathematics, Yuquan Campus, Zhejiang
University, Hangzhou 310027, China \newline E-mail address:
$^1$lxzhang@mail.hz.zj.cn  \quad $^2$ hwmath2003@hotmail.com
 }}
\author{Li-Xin Zhang$^1$ and Wei Huang$^2$}

\bigskip


\maketitle

\bigskip

\bigskip


\bigskip
\begin{center}
{\bf Abstract}
\end{center}

\noindent In literature, the central limit theorems  for the product
of sums of various random variables have studied. The purpose of
this note is to show that this kind of results are corollary of the
invariance principle.



\bigskip
\bigskip
\setcounter{equation}{0}

Let $\{X_k;k\ge 1\}$ be a sequence of i.i.d  exponential random
variables with mean $1$, $S_n=\sum_{k=1}^nX_k$, $n \ge 1$. Arnold
and Villase\~{n}or (1998) proved that
\begin{equation}\label{eq1}
\left(\prod_{k=1}^n \frac{S_k}{k}\right)^{1/\sqrt{n}}\overset{\mathscr D}\to
e^{\sqrt{2} N(0,1)}, \ \ \ \text{as}  \ \ n \to \infty,
\end{equation}
where $N(0,1)$ is a standard normal random variable. Later Rempala
and Wesolowski (2002) extended such a central limit theorem to
general i.i.d. positive random variables. Recently, the central
limit theorems for product of sums  have been also studied for
dependent random variables (c.f., Gonchigdanzan and Rempala (2006)).
In this note, we will show that this kind of results can follow from
the invariance principle.

Let $\{S_n; n\ge 1\}$ be a sequence of positive random variables. To
present our main idea, we assume that (possibly in an enlarged probability space in which the sequence
$\{S_n; n\ge 1\}$ is redefined
without changing its distribution) there exists a standard Wiener
process $\{W(t): t \ge 0\}$ and two positive constants $\mu$ and
$\sigma$ such that
\begin{equation}\label{eq2}
S_n-n\mu -\sigma W(n)=o (\sqrt{n}) \;\;\; a.s.
\end{equation}
Then
\begin{align}\label{eqapp}
\log \prod_{k=1}^n \frac{S_k}{k\mu}
=&\sum_{k=1}^n \log \frac{S_k}{k\mu}
=\sum_{k=1}^n \log\left(1+\frac{\sigma}{\mu}\frac{W(k)}{k}+o(k^{-1/2})\right)
\nonumber\\
=& \sum_{k=1}^n \left(\frac{\sigma}{\mu}\frac{W(k)}{k}+o(k^{-1/2})\right)=
\frac{\sigma}{\mu}\sum_{k=1}^n \frac{W(k)}{k}+ o(\sqrt{n})
\nonumber \\
=& \frac{\sigma}{\mu}\int_0^n \frac{W(x)}{x} dx+ o(\sqrt{n})\;\;\;
a.s.,
\end{align}
where $\log x= ln (x\vee e)$.  It follows that
$$ \frac{\mu}{\sigma}\frac{1}{\sqrt{n}}\log \prod_{k=1}^n \frac{S_k}{k\mu}
\overset{\mathscr D}\to \int_0^1\frac{W(x)}{x}dx, \;\;\;\text{as} \;\; n \to
\infty.$$
It is easy seen that the
random variable in the right hand side is a normal random variable with
$$ \ep \int_0^1\frac{W(x)}{x}dx=\int_0^1\frac{\ep W(x)}{x}dx=0 $$
and
$$ \ep\left(\int_0^1\frac{W(x)}{x}dx\right)^2=\int_0^1\int_0^1\frac{\ep W(x) W(y)}{xy}dxdy
=\int_0^1\int_0^1\frac{\min(x,y)}{xy}dxdy=2.$$
So
\begin{equation}\label{eq3}
\left(\prod_{k=1}^n
\frac{S_k}{k\mu}\right)^{\gamma/\sqrt{n}}\overset{\mathscr D}\to e^{\sqrt{2}
N(0,1)},\ \ \ \text{as} \ \ n \to \infty,
\end{equation}
where $\gamma =\mu/\sigma$. If $S_n$ is the partial  sum of a
sequence $\{X_k; k\ge 1\}$ of i.i.d. random variables, then
(\ref{eq2}) is satisfied when $\ep |X_k|^2\log\log |X_k|<\infty$.
(\ref{eq2}) is known as the strong invariance principle. To show
(\ref{eq3}) holds for sums of  i.i.d. random variables only with the
finite second moments, we replace the condition (\ref{eq2}) by a
weaker one. The following is our main result.

\begin{theorem}\label{thm1.1}
Let $\{S_k; k\ge 1\}$ be a nondecreasing sequence of positive random
variables. Suppose there exists a standard Wiener process $\{W(t); t
\ge 0\}$ and two positive constants $\mu$ and $\sigma$ such that
\begin{equation}\label{eqth1.1}
W_n(t)=:\frac{S_{[nt]}-[nt]\mu}{\sigma\sqrt{n}}\overset{\mathscr D}\to W(t)
\;\text{ in }\; D[0,1], \;\; \text{as} \;\; n \to \infty
\end{equation}
 and
\begin{equation}\label{eqth1.2}
\sup_n \frac{E|S_n-n\mu |}{\sqrt{n}}<\infty.
\end{equation}
Then
\begin{equation}\label{eqth1.3}
\left(\prod_{k=1}^{[nt]} \frac{S_k}{k\mu}\right)^{\gamma/\sqrt{n}}
\overset{\mathscr D}\to \exp\left\{\int_0^t \frac{W(x)}{x}dx\right\} \;\text{
in }\; D[0,1], \;\; \text{as} \;\; n \to \infty,
\end{equation}
where $\gamma =\mu/\sigma$.
\end{theorem}

\begin{remark} (\ref{eqth1.1}) is known as the weak invariance principle.
The conditions (\ref{eqth1.1}) and (\ref{eqth1.2}) are satisfied for
many random variables sequences. For example, if $\{X_k;k\ge 1\}$
are i.i.d. positive random variables with mean $\mu$ and  variance
$\sigma^2$ and
 $S_n=\sum_{i=1}^n X_k$,
then (\ref{eqth1.1}) is satisfied by the invariance principle (c.f., Theorem 14.1 of Billingsley (1999)). Also,
$$ \ep\left[\frac{|S_n-n\mu|}{\sqrt{n}}\right]\le \left\{\Var\left[\frac{S_n-n\mu}{\sqrt{n}}\right]\right\}^{1/2}
=\sigma, $$ by Cauchy-Schwarz inequality. Condition (\ref{eqth1.2})
is also satisfied. Many dependent random sequences also satisfy
these two conditions.
\end{remark}

{\bf Proof of Theorem \ref{thm1.1}.} For $x>-1$, write
$\log(1+x)=x+x\theta(x)$, where $\theta(x)\to 0$, as $x\to 0$. Then
for any $t>0$,
\begin{align}\label{eqproofth1.1}
&\log\left(\prod_{k=1}^{[nt]} \frac{S_k}{k\mu}\right)^{\gamma/\sqrt{n}}
= \frac{1}{\sigma\sqrt{n}}\sum_{k=1}^{[nt]} \frac{S_k-k\mu}{k}
+ \frac{1}{\sigma\sqrt{n}}\sum_{k=1}^{[nt]} \frac{S_k-k\mu}{k}\theta\left(\frac{S_k}{k\mu}-1\right).
\end{align}
Notice that for any $\rho>1$,
\begin{align*}
\max_{\rho^n\le k\le \rho^{n+1}}\frac{|S_k-k\mu|}{k}\le \max\left\{\frac{|S_{\rho^{n+1}}-\rho^{n+1}\mu|}{\rho^n},
\frac{|S_{\rho^n}-\rho^n\mu|}{\rho^n}\right\}+(\rho-1)\mu.
\end{align*}
Together with (\ref{eqth1.2}), it follows that, for any $n_0 \ge 1$,
\begin{align*}
&\ep \left[\max_{k\ge \rho^{n_0}}\frac{|S_k-k\mu|}{k}\right]
\le \rho\ep \left[\max_{n\ge n_0} \frac{|S_{\rho^n}-\rho^n\mu|}{\rho^n}\right]+(\rho-1)\mu\\
\le &
\rho\sup_k\frac{E|S_k-k\mu|}{\sqrt{k}}\sum_{n=n_0}^{\infty}\rho^{-n/2}+(\rho-1)\mu\to
0,
\end{align*}
as $n_0\to \infty$ and then $\rho\to 1$. It follows that
$$\max_{k\ge k_0}\left|\frac{S_k}{k\mu}-1\right|\overset{P}\to 0, \;\;\text{ as }\;\; k_0\to \infty,$$
which implies that
$$ \frac{S_k}{k\mu}-1 \to 0 \;\; a.s., \ \ \  \text{as} \ \ k\to \infty. $$
Hence we conclude that
$$\theta\left(\frac{S_k}{k\mu}-1\right)\to 0 \;\; a.s., \ \ \  \text{as} \ \ k\to \infty. $$
On the other hand, by (\ref{eqth1.2}), we have
\begin{equation}\label{eqproofth1.2}
 \frac{1}{\sqrt{n}}\ep\left[\sum_{k=1}^n \frac{|S_k-k\mu|}{k}\right]
\le C_0 \frac{1}{\sqrt{n}}\sum_{k=1}^n \frac{1}{\sqrt{k}}\le 2C_0.
\end{equation} It follows that
$$\max_{0 \le t\le 1} \left|\frac{1}{\sigma\sqrt{n}}\sum_{k=1}^{[nt]} \frac{S_k-k\mu}{k\mu}\theta\left(\frac{S_k}{k\mu}-1\right)\right|
= \frac{1}{\sqrt{n}}\sum_{k=1}^n \frac{|S_k-k\mu|}{k}o(1)=o_P(1). $$
So, according to (\ref{eqproofth1.1}) it is suffices to show that
\begin{equation}\label{eqproofth1.3}
Y_n(t)=:\frac{1}{\sigma\sqrt{n}}\sum_{k=1}^{[nt]}
\frac{S_k-k\mu}{k}\overset{\mathscr D}\to \int_0^t \frac{W(x)}{x}dx \;\text{
in } \; D[0,1], \ \ \ \text{as} \ \ n \to \infty.
\end{equation}
Let
$$ H_{\epsilon}(f)(t)=\begin{cases}\displaystyle
  \int_{\epsilon}^t \frac{f(x)}{x}, & t >\epsilon,\\
  0, &  0\le t\le \epsilon
  \end{cases}
  $$
  and
$$ Y_{n,\epsilon}(t)=\begin{cases}\displaystyle \frac{1}{\sigma\sqrt{n}}\sum_{k=[n\epsilon]+1}^{[nt]}
\frac{S_k-k\mu}{k}, & t>\epsilon, \\
0, & 0\le t\le \epsilon. \end{cases} $$
It is obvious that
\begin{equation}\label{eqproofth1.4}
\max_{0\le t\le 1}\left|\int_0^t \frac{W(x)}{x}dx-H_{\epsilon}(W)(t)\right|=\sup_{0\le t\le \epsilon} \left| \int_0^t \frac{W(x)}{x}dx\right|\to
0 \quad a.s., \text{ as } \epsilon \to 0
\end{equation}
and
\begin{align}\label{eqproofth1.5} &
E\max_{0\le t\le \epsilon}\left|Y_n(t)-Y_{n,\epsilon}(t)\right|
=E\bigg\{\max_{0\le t\le \epsilon}E\Big|\frac{1}{\sigma\sqrt{n}}\sum_{k=1}^{[nt]}
 \frac{S_k-k\mu}{k}\Big|\bigg\}
\nonumber \\
\le & \frac{1}{\sigma\sqrt{n}}\sum_{k=1}^{[n\epsilon]}
\frac{\ep|S_k-k\mu|}{k}
 \le
\frac{C_0}{\sigma\sqrt{n}}\sum_{k=1}^{[n\epsilon]} \frac{1}{\sqrt k}
 \le \frac{2C_0}{\sigma\sqrt{n}}
 \sqrt{[n\epsilon]} \le C \sqrt{\epsilon},
\end{align}
by (\ref{eqth1.2}). On the other hand, it is easily seen that
\begin{align*}
&\sup_{\epsilon\le t\le 1}\bigg| \sum_{k=[n\epsilon]+1}^{[nt]} \frac{S_k-k\mu}{k}-
\int_{n\epsilon}^{nt}\frac{S_{[x]}-[x]\mu}{x}dx\bigg|
\\
=&\sup_{\epsilon\le t\le 1}\bigg| \int_{[n\epsilon]+1\le x<[nt]+1}
\frac{S_{[x]}-[x]\mu}{[x]}dx- \int_{n\epsilon}^{nt}\frac{S_{[x]}-[x]\mu}{x}dx\bigg|\\
\le &\bigg| \int_{n\epsilon \le x<[n\epsilon]+1}
\frac{S_{[x]}-[x]\mu}{x}dx\bigg|
+\sup_{\epsilon\le t\le 1}\bigg| \int_{nt
\le x<[nt]+1} \frac{S_{[x]}-[x]\mu}{x}dx\bigg|\\
&+\sup_{\epsilon\le t\le 1}\bigg|\int_{[n\epsilon]+1 \le x<[nt]+1}\big(S_{[x]}-[x]\mu\big)\bigg(\frac{1}{x}-\frac{1}{[x]}\bigg)dx\bigg|\\
\le & \max_{k\le n} |S_{k}-k\mu|\sup_{\epsilon\le t\le 1}\bigg(
\frac{2}{n\epsilon}+\frac{2}{n t}+\frac{1}{n\epsilon}\bigg)\\
\le &
5\max_{k\le n} |S_{k}-k\mu|/(n\epsilon)=O_P(\sqrt{n})/n=o_P(1)
\end{align*}
by (\ref{eqth1.1}). So
\begin{align*}\frac{1}{\sigma\sqrt{n}}\sum_{k=[n\epsilon]+1}^{[nt]} \frac{S_k-k\mu}{k}
=&\frac{1}{\sigma\sqrt{n}}\int_{n\epsilon}^{nt}\frac{S_{[x]}-[x]\mu}{x}dx+o_P(1)
=\int_{\epsilon}^{t}\frac{W_{[nx]}}{x}dx+o_P(1)
\end{align*}
uniformly in $t\in[\epsilon, 1]$.
Notice that $H_{\epsilon}(\cdot)$ is a continuous mapping on the space $D[0,1]$. Using the continuous mapping
theorem (c.f., Theorem 2.7 of Billingsley (1999)) it follows that
\begin{equation}\label{eqproofth1.6}
Y_{n,\epsilon}(t)=H_{\epsilon}(W_n)(t)+o_P(1)\overset{\mathscr D}\to H_{\epsilon}(W)(t) \;\text{ in
} \; D[0,1], \;\;\;\text{as} \;\; n \to
\infty.
\end{equation}
Combining (\ref{eqproofth1.4})--(\ref{eqproofth1.6}) yields
(\ref{eqproofth1.3}) by Theorem 3.2 of Billingsley (1999). $\Box$

\begin{theorem}\label{thm1.2}
Let $\{S_k; k\ge 1\}$ be a sequence of positive random variables.
Suppose there exists a standard Wiener process $\{W(t); t \ge 0\}$
and two positive constants $\mu$ and $\sigma$ such that
\begin{equation}\label{eqth2.1}
S_n-n\mu-\sigma W(n) = o\big(\sqrt{n\log\log n}\big)\;\; a.s.
\end{equation}
Let
$$\Cal F=\left\{f(t)=\int_0^tf^{\prime}(u)du: f(0)=0,
\int_0^1(f^{\prime}(u))^2 du\le  1, 0\le u\le 1\right\} $$
be the set of  continuous functions on $[0,1]$ which are absolutely continuous  (with
respect to Lebesgue measure) and have derivatives in the unit ball of $L_2[0,1]$.
 Then with probability one
\begin{equation}\label{eqth2.2}
\left\{\Big(\prod_{k=1}^{[nt]}
\frac{S_k}{k\mu}\Big)^{\gamma/\sqrt{2n\log\log n}}; 0\le t\le
1\right\}_{n=3}^{\infty} \text{ is relatively compact }
\end{equation}
and the limit set is
$$ \left\{ \exp\Big\{\int_0^x \frac{f(u)}{u} du\Big\}: f\in \Cal F, 0\le x\le 1\right\}. $$
 Particularly,
\begin{equation}\label{eqth2.3}
\limsup_{n\to \infty} \left(\prod_{k=1}^n
\frac{S_k}{k\mu}\right)^{\gamma/\sqrt{2n\log\log
n}}=e^{\sqrt{2}}\;\; a.s.
\end{equation}
\end{theorem}

{\bf Proof of Theorem \ref{thm1.2}.} Similar to (\ref{eqapp}), we
have
$$\log \prod_{k=1}^n \frac{S_k}{k\mu}
= \frac{\sigma}{\mu}\int_0^n \frac{W(x)}{x} dx+ o(\sqrt{n\log\log n})\;\; a.s. $$
Notice
$$ \frac{1}{\sqrt{2n\log\log n}}\int_0^{nt} \frac{W(x)}{x} dx =\int_0^t \frac{1}{u}\frac{W(nu)}{\sqrt{2n\log\log n}}du $$
and with probability one
$$\left\{ \frac{W(nt)}{\sqrt{2n\log\log n}}: 0\le t\le 1\right\}_{n=3}^{\infty} \; \text{ is relatively compact }$$
with $\Cal F$ being the limit set (c.f., Theorem 1.3.2 of Cs\H org\"
o and R\'ev\'esz (1981) or Strassen (1964)). The first part of
conclusion follows immediately. For (\ref{eqth2.3}), it suffices to
show that
\begin{equation}\label{eqproofth2.1}
 \sup_{f\in \Cal F} \sup_{0\le t\le 1} \int_0^t \frac{f(u)}{u}du \le \sqrt{2}
 \end{equation}
and
\begin{equation}\label{eqproofth2.2}
 \sup_{f\in \Cal F} \int_0^1 \frac{f(u)}{u}du \ge \sqrt{2}.
 \end{equation}
For any $f\in \Cal F$, using Cauchy-Schwarz inequality, we have
\begin{align*}
\int_0^t \frac{f(u)}{u}du =&\int_0^t \frac{1}{u} \int_0^u f^{\prime} (v) dvdu
=\int_0^t \int_{v}^t f^{\prime} (v)\frac{1}{u}   dudv \\
=&\int_0^t f^{\prime} (v)\log\frac{t}{v} dv \le \left( \int_0^t \left(\log\frac{t}{v}\right)^2 dv\right)^{1/2}
\left(\int_0^t \big(f^{\prime} (v) \big)^2 dv\right)^{1/2} \\
\le & \left( \int_0^t \left(\log\frac{t}{v}\right)^2
dv\right)^{1/2}=\sqrt{2t }\le \sqrt{2},
\end{align*}
where $0\le t\le 1$. Then (\ref{eqproofth2.1}) is proved. Now, let
$f(t)=(t-t\log t)/\sqrt{2}$, $f(0)=0$. Then $ f\in \Cal F$ and
$$ \int_0^1\frac{f(u)}{u} du = \frac{1}{\sqrt{2}}\int_0^1( 1-\log u)du =\sqrt{2}. $$
(\ref{eqproofth2.2}) is proved. $\Box$.
\bigskip


\begin{thebibliography}{99}

\bibitem {}  B.C. Arnold and J. A. Villase\~{n}or. The Asymptotic Distributions of Sums of Records.
{\em Extremes}, {\bf 1}, {No.3} (1998), 351-363. MR1814709
(2002a:60025)

\bibitem {} P. Billingsley. {\em Convergence of Probability Measures}, Joh Wiley \&
Sons, INC, New York (1999).

\bibitem{} M. Cs\"org\H o and P. R\'ev\'esz. {\em Strong
Approximations in Probability and Statistics}, Akad\'{e}miai
Kiad\'{o}, Budapest (1981).

\bibitem {} K. Gonchigdanzan and G. A. Rempala. A note on the almost sure limit theorem
for the product of partial sums. {\em Applied Math. Lett.}, {\bf
19}, {No. 2} (2006), 191-196. MR2198407

\bibitem {} G. Rempala and J. Wesolowski. Asymptotics for products of sums and
U-statistics. {\em Elect. Comm. in Prob.}, {\bf 7} (2002), 47-54.
MR1887173 (2002k:60070)


\bibitem {} V. Strassen. An invariance principle for the law
of the iterated logarithm. {\em Z. Wahrsch. Keitsth. verw. Gebiete},
{\bf 3} (1964), 211-226. MR0175194 (30 \#5379)
\end{thebibliography}
\end{document}